\documentclass[a4paper,11pt]{article}
\usepackage{amssymb}

\newcommand{\R}{\mathbb{R}}

\newcommand{\N}{\mathbb{N}}

\newcommand{\beq}{\begin{equation} }
\newcommand{\eqq}{\end{equation} }
\newcommand{\cuad}{{\sqcap\kern-.68em\sqcup}}

\newcommand{\norm}[1]{\|#1\|}

\newtheorem{definition}{Definition}[section]
\newtheorem{teo}{Theorem}[section]

\newtheorem{proposition}{Proposition}[section]

\newtheorem{lemma}{Lemma}[section]

\newtheorem{remark}{Remark}[section]
\newcommand{\bremark}{\begin{remark} \em}
\newcommand{\eremark}{\end{remark} }

\begin{document}

\begin{center}{\bf  \Large    Elliptic equations involving
general subcritical

\medskip
 source nonlinearity
  and  measures
 }\bigskip
%%%%%%%%%%%%%%%%%%%%%%%%%%%%%%%%%%%%%%%%%%%%%%%%%%%%%%%%%%%%%%%%%%%%%%
%%%%%%%%%%%%%%%%%%%%%%%%%%%%%%%%%%%%%%%%%%%%%%%%%%%%%%%%%%%%%%%%%%%%%%

 Huyuan Chen\footnote{hc64@nyu.edu} \qquad Patricio Felmer\footnote{pfelmer@dim.uchile.cl}
\qquad Laurent V\'{e}ron\footnote{Laurent.Veron@lmpt.univ-tours.fr}
\bigskip

\medskip

\begin{abstract}
In this article, we study the existence of positive solutions to  elliptic equation (E1)
$$(-\Delta)^\alpha u=g(u)+\sigma\nu \quad{\rm in}\quad \Omega,$$
subject to the condition (E2)
$$u=\varrho\mu\quad {\rm on}\quad \partial\Omega\ \ {\rm if}\ \alpha=1\qquad {\rm or\ \ in}\ \ \Omega^c \ \ {\rm if}\ \alpha\in(0,1),$$
 where $\sigma,\varrho\ge0$,  $\Omega$ is an open bounded $C^2$ domain in $\R^N$,  $(-\Delta)^\alpha$
denotes the fractional Laplacian with $\alpha\in(0,1)$ or Laplacian operator if $\alpha=1$, $\nu,\mu$ are suitable Radon measures
 and $g:\R_+\mapsto\R_+$ is a continuous  function.

 We introduce an approach to obtain  weak solutions
for problem (E1)-(E2) when  $g$ is integral subcritical and $\sigma,\varrho\ge0$ small enough.
\end{abstract}
\end{center}

%\tableofcontents \vspace{1mm}
  \vspace{1mm}


\setcounter{equation}{0}
\section{Introduction}
Let $\alpha\in(0,1]$, $\Omega$ be an open bounded $C^2$ domain in $\R^N$ with $N>2\alpha$,  $\rho(x)=dist(x,\partial\Omega)$, $g:\R_+\mapsto\R_+$ be a continuous  function  and denote by
$(-\Delta)^\alpha$ the Laplacian operator if $\alpha=1$ or  the fractional Laplacian with $\alpha\in(0,1)$ defined as
$$(-\Delta)^\alpha  u(x)=\lim_{\varepsilon\to0^+} (-\Delta)_\varepsilon^\alpha u(x),$$
where for $\varepsilon>0$,
$$
(-\Delta)_\varepsilon^\alpha  u(x)=-\int_{\R^N}\frac{ u(z)-
u(x)}{|z-x|^{N+2\alpha}}\chi_\varepsilon(|x-z|) dz
$$
and
$$\chi_\varepsilon(t)=\left\{ \arraycolsep=1pt
\begin{array}{lll}
0,\quad & {\rm if}\quad t\in[0,\varepsilon],\\[2mm]
1,\quad & {\rm if}\quad t>\varepsilon.
\end{array}
\right.$$

 Our first purpose of this paper is to study the
existence of weak solutions  to the semilinear  elliptic problem
\begin{equation}\label{eq1.1}
 (-\Delta)^\alpha  u=g(u)+\sigma \nu\quad  {\rm in}\quad\Omega,
\end{equation}
subject to the Dirichlet boundary condition
\begin{equation}\label{eq1.01}
u=0\quad {\rm on}\quad \partial\Omega\ \ {\rm if}\ \alpha=1\quad {\rm or\  in}\ \ \Omega^c \ \ {\rm if}\ \alpha\in(0,1),
\end{equation}
where  $\sigma>0$,  $\nu\in\mathfrak{M}(\Omega,\rho^\beta)$ with $\beta\in[0,\alpha]$ and $\mathfrak{M}(\Omega,\rho^\beta)$ being the
space of Radon measures in $\Omega$ satisfying
$$
\int_{\Omega}\rho^\beta d|\nu|<+\infty.
$$
In particular, we denote
$\mathfrak{M}^b(\Omega)=\mathfrak{M}(\Omega,\rho^0)$. The associated positive cones are respectively
$\mathfrak{M}_+(\Omega,\rho^\beta)$ and $\mathfrak{M}_+^b(\Omega)$.

When $\alpha=1$, problem (\ref{eq1.1})-(\ref{eq1.01})  has been studied for some decades. The basic method developed by Ni \cite{N} and
Ratto-Rigoli-V\'{e}ron \cite{RRV} is
to iterate
$$u_{n+1}=\mathbb{G}_1[g(u_n)]+\sigma\mathbb{G}_1[\nu],\quad \forall n\in\N.$$
The crucial ingredient in this approach is to derive a function $v$ satisfying
$$v\ge \mathbb{G}_1[g(v)]+\sigma\mathbb{G}_1[\nu].$$
Later on,  Baras-Pierre \cite{BP} applied duality argument  to derive weak solution of problem (\ref{eq1.1})-(\ref{eq1.01}) with $\alpha=1$ under the hypotheses:\\
$(i)$
the mapping $r\mapsto g(r)$ is  nondecreasing,  convex  and continuous;\\
$(ii)$
there exist $c_0>0$ and $\xi_0\in C_0^{1.1}(\Omega)$, $\xi_0\not=0$ such that
    $$g^*\left(c_0\frac{-\Delta \xi_0}{\xi_0}\right)\in L^1(\Omega),$$
    where $g^*$ is the conjugate function of $g$;\\
$(iii)$
$$
\int_\Omega \xi d\nu\le \int_\Omega g^*\left(\frac{-\Delta \xi}{\xi}\right)dx, \qquad \forall \xi\in C_0^{1.1}(\Omega).
$$
When  $g$ is pure power source, Brezis-Cabr\'{e} \cite{BC} and Kalton-Verbitsky \cite{KV} pointed out that
the necessary condition for existence of weak solution to
\begin{equation}\label{eq1.2}
 \arraycolsep=1pt
\begin{array}{lll}
 -\Delta   u=u^p+\sigma \nu\quad & {\rm in}\quad\Omega,\\[2mm]
 \phantom{   -\Delta  }
u=0\quad & {\rm on}\quad \partial\Omega,
\end{array}
\end{equation} is that
\begin{equation}\label{eq1.3}
\mathbb{G}_1[(\mathbb{G}_1[ \nu])^p]\le c_1\mathbb{G}_1[ \nu],
\end{equation}
for some  $c_1>0$.
Bidaut-V\'{e}ron and Vibier in \cite{BV} proved that (\ref{eq1.3}) holds for $p<\frac{N+\beta}{N+\beta-2}$ and problem
(\ref{eq1.2})
admits a weak solution if $\sigma>0$ small.
While it is not easy to get explicit condition for general nonlinearity by  above methods.

In this article, we introduce a new method to obtain the weak solution of problem (\ref{eq1.1})-(\ref{eq1.01}) involving general nonlinearity without convex and nondecreasing properties, which is inspired by the Marcinkiewicz spaces approach.

Let us first make precise the definition of weak solution to (\ref{eq1.1})-(\ref{eq1.01}).
\begin{definition}\label{weak definition}
We say that $u$ is a weak solution of (\ref{eq1.1})-(\ref{eq1.01}), if $u\in
L^1(\Omega)$,  $g(u)\in L^1(\Omega,\rho^\alpha dx)$  and
$$
\int_\Omega u(-\Delta)^\alpha\xi dx=\int_\Omega g(u)\xi dx+\sigma\int_\Omega\xi
d\nu,\quad\  \forall \xi\in \mathbb{X}_{\alpha},
$$
where $\mathbb{X}_{\alpha}=C_0^{1.1}(\Omega)$ if $\alpha=1$ or $\mathbb{X}_{\alpha}\subset C(\R^N)$ with $\alpha\in(0,1)$ is the space of functions
$\xi$ satisfying:\smallskip

\noindent (i) ${\rm supp}(\xi)\subset\bar\Omega$,\smallskip

\noindent(ii) $(-\Delta)^\alpha\xi(x)$ exists for all $x\in \Omega$
and $|(-\Delta)^\alpha\xi(x)|\leq C$ for some $C>0$,\smallskip

\noindent(iii) there exist $\varphi\in L^1(\Omega,\rho^\alpha dx)$
and $\varepsilon_0>0$ such that $|(-\Delta)_\varepsilon^\alpha\xi|\le
\varphi$ a.e. in $\Omega$, for all
$\varepsilon\in(0,\varepsilon_0]$.\smallskip
\end{definition}

We denote by  $G_\alpha$ the Green kernel of $(-\Delta)^\alpha$ in
$\Omega\times\Omega $ and  by $\mathbb{G}_\alpha[.]$ the associated Green operator
defined by
$$
\mathbb{G}_\alpha[\nu](x)=\int_{\Omega}G_\alpha(x,y) d\nu(y),\qquad\forall  \nu\in
\mathfrak{M}(\Omega,\rho^\alpha).
$$
Our first result states as follows.

\begin{teo}\label{teo 1}
Let $\alpha\in(0,1]$, $\sigma>0$ and $\nu\in\mathfrak{M}_+(\Omega,\rho^{\beta})$ with $\beta\in[0,\alpha]$.

$(i)$ Suppose that
\begin{equation}\label{06-08-2}
g(s)\le c_2s^{p_0}+\epsilon,\quad \forall s\ge0,
\end{equation}
for some $p_0\in(0,1]$, $c_2>0$ and $\epsilon>0$.
Assume more that  $c_2$ is small enough when $p_0=1$.

Then problem (\ref{eq1.1})-(\ref{eq1.01}) admits a  weak nonnegative solution $u_\nu$ which satisfies
\begin{equation}\label{1.5}
 u_\nu\ge \sigma\mathbb{G}_\alpha[\nu].
\end{equation}

$(ii)$ Suppose that
\begin{equation}\label{1.1}
 g(s)\le c_3s^{p_*}+\epsilon,\quad \forall s\in[0,1]
\end{equation}
and
\begin{equation}\label{1.4}
g_\infty:=\int_1^{+\infty} g(s)s^{-1-p^*_{\beta}}ds<+\infty,
\end{equation}
where  $c_3,\epsilon>0$, $p_*>1$ and
$
p^*_\beta= \frac{N+\beta}{N-2\alpha+\beta}.
$

Then there exist $\sigma_0,\epsilon_0>0$  depending on $c_3, p_*, g_\infty$ and $ p_\beta^*$  such that for $\sigma\in[0,\sigma_0)$ and $\epsilon\in(0,\epsilon_0)$, problem
(\ref{eq1.1})-(\ref{eq1.01}) admits a nonnegative  weak solution $u_\nu$ which satisfies (\ref{1.5}).
\end{teo}

We remark that $(i)$ we do not require any restriction on parameters $c_2, \epsilon, \sigma$ when $p_0\in(0,1)$ or
on parameters  $ \epsilon, \sigma$ when $p_0=1$;
$(ii)$ the assumption (\ref{1.4}) is called as integral subcritical condition, which is
usually used in dealing with elliptic problem with absorption nonlinearity and measures, see the references \cite{BV,CV1,CV2,V}.

Let us sketch  the proof  of Theorem \ref{teo 1}. We first approximate the nonlinearity $g$ and Radon measure
$\nu$ by $\{g_n\}$ and $\{\nu_n\}$ respectively, then we make use of   the Marcinkiewicz properties and embedding theorems
 to obtain  that for $n\geq 1$, problem
$$
 (-\Delta)^\alpha u_{n}= g_{n}(u_{n})+\sigma\nu_n\quad  {\rm in}\quad\Omega,
$$
subject to condition (\ref{eq1.01}),
admits a nonnegative solution $u_n$ by Schauder's fixed point theorem. The crucial point is to obtain uniformly bound of $\{u_n\}$ in  the Marcinkiewicz space. The proof ends by getting a subsequence of $\{u_n\}$ that converges in the sense of Definition \ref{weak definition}.

Our second purpose in this note is to obtain the weak solution to elliptic equations involving boundary measures.  Firstly, we study the weak solution of  %It is known
%that the fractional equation involving boundary value is ill-posed.
\begin{equation}\label{eq1.1b} \arraycolsep=1pt
\begin{array}{lll}
-\Delta u=g(u) \quad & {\rm in}\quad \Omega,
\\[2mm]\phantom{--}
u=\varrho \mu\quad &{\rm on}\quad \partial\Omega,
\end{array}
\end{equation}
where  $\varrho>0$ and $\mu\in \mathfrak{M}^b_+(\partial\Omega)$ the space of nonnegative bounded Radon measure on $\partial\Omega$.
When $g(s)=s^p$ with $p<\frac{N+1}{N-1}$, the weak solution to problem (\ref{eq1.1b}) is derived by Bidaut-V\'{e}ron and Vivier in \cite{BV}
by using iterating procedure. More interests on boundary measures refer to \cite{BM,BY,GV,MV1,MV2,MV}.

\begin{definition}\label{weak definition b}
We say that $u$ is a weak solution of (\ref{eq1.1b}) , if $u\in
L^1(\Omega)$,  $g(u)\in L^1(\Omega,\rho dx)$  and
$$
\int_\Omega u(-\Delta)\xi dx=\int_\Omega g(u)\xi dx+\varrho\int_{\partial\Omega}\frac{\partial\xi(x)}{\partial \vec n_x}
d\mu(x),\quad  \forall \xi\in C_0^{1.1}(\Omega),
$$
where $\vec n_x$ is the unit normal vector pointing outside of $\Omega$ at point $x$.
\end{definition}

We denote by  $P$ the Poisson kernel of $-\Delta$ in
$\Omega\times\partial\Omega $ and  by $\mathbb{P}[.]$ the associated Poisson operator
defined by
$$
\mathbb{P}[\mu](x)=\int_{\partial\Omega}P(x,y) d\mu(y),\qquad\forall \mu\in
\mathfrak{M}^b(\partial\Omega).
$$
Our second result states as follows.

\begin{teo}\label{teo 2}
Let $\varrho>0$ and $\mu\in\mathfrak{M}_+^b(\partial\Omega)$.

$(i)$ Suppose that
\begin{equation}\label{06-08-2}
g(s)\le c_4s^{q_0}+\epsilon,\quad \forall s\ge0,
\end{equation}
for some $q_0\in(0,1]$, $c_4>0$ and $\epsilon>0$.
Assume more that  $c_4$ is small enough when $q_0=1$.

Then problem (\ref{eq1.1b}) admits a  weak nonnegative solution $u_\mu$ which satisfies
\begin{equation}\label{1.6}
 u_\mu\ge \varrho\mathbb{P}[\mu].
\end{equation}

$(ii)$ Suppose that
\begin{equation}\label{g10}
 g(s)\le c_5s^{q_*}+\epsilon,\quad \forall s\in[0,1]
\end{equation}
and
\begin{equation}\label{g1+}
g_{\infty}:=\int_1^{+\infty} g(s)s^{-1-q^*}ds<+\infty,
\end{equation}
where $c_5,\epsilon>0$, $q_*>1$ and $
q^*= \frac{N+1}{N-1}.
$

Then there exist $\varrho_0,\epsilon_0>0$  depending on $c_5, q_*, g_\infty$ and $ q^*$  such that for $\varrho\in[0,\varrho_0)$ and $\epsilon\in[0,\epsilon_0)$, problem
(\ref{eq1.1b}) admits a nonnegative  weak solution $u_\mu$ which satisfies
(\ref{1.6}).

\end{teo}

We remark that the key-point in the proof of Theorem \ref{teo 2}
is to derive the uniform bound in Marcinkiewicz quasi-norm to the solutions of
\begin{equation}\label{1.7}
\arraycolsep=1pt
\begin{array}{lll}
 -\Delta  u= g_n(u+\varrho\mathbb{P}[\mu])\quad & {\rm in}\quad\Omega,\\[2mm]
 \phantom{   -\Delta  }
u=0\quad & {\rm on}\quad \partial\Omega,
\end{array}
\end{equation}
where $\{g_n\}$ is a sequence of $C^1$ bounded functions approaching to $g$ in $L^\infty_{loc}(\R_+)$.
In fact, the weak solution $u_\mu$ could be decomposed into
$$u_\mu=v_\mu+\varrho \mathbb{P}[\mu],$$
where $v_\mu$ is a weak solution to (\ref{1.7}) replaced $g_n$ by $g$.

Inspired by the fact above, we  give the definition of  weak solution  to
\begin{equation}\label{eq1.1bt}
\arraycolsep=1pt
\begin{array}{lll}
(-\Delta)^\alpha u=g(u) \quad & {\rm in}\quad \Omega,
\\[2mm]\phantom{---}
u=\varrho\mu\quad &{\rm in}\quad \Omega^c
\end{array}
\end{equation}
as follows.

\begin{definition}\label{weak definition bt}
We say that $u_\mu$ is a weak solution of (\ref{eq1.1bt}) ,
 if
 $$u_\mu=v_\mu+\varrho \mathbb{G}_\alpha[ w_\mu],$$
where
\begin{equation}\label{poisson type}
w_\mu(x) =\int_{\Omega^c}\frac{d\mu(z)}{|z-x|^{N+2\alpha}},\quad x\in\Omega
\end{equation}
and $v_\mu$ is a solution of
 \begin{equation}\label{eq1.1btc}
\arraycolsep=1pt
\begin{array}{lll}
(-\Delta)^\alpha u=g(u+\varrho \mathbb{G}_\alpha[ w_\mu]) \quad & {\rm in}\quad \Omega,
\\[2mm]\phantom{---}
u=0 \quad &{\rm in}\quad \Omega^c
\end{array}
\end{equation}
in the  sense of Definition \ref{weak definition}.
\end{definition}
In Definition \ref{weak definition bt}, the function $\mathbb{G}_\alpha[ w_\mu]$ plays the role of
$\mathbb{P}[\mu]$ when $\alpha=1$.
In order to better classify the measures  tackled in follows,
 we denote
\begin{equation}\label{5.2}
\mathfrak{R}_\beta:=\{\mu\in \mathfrak{M}_+(\Omega^c):\  w_\mu\in L^1(\Omega,\rho^\beta dx)\},
\end{equation}
where $ \beta\in[0,\alpha]$  and $w_\mu$ is given by (\ref{poisson type}).

\begin{teo}\label{teo 3}

Let $\alpha\in(0,1)$, $\sigma>0$ and $\mu\in \mathfrak{R}_{\beta}$ with $ \beta\in[0,\alpha]$.

$(i)$ Suppose that
\begin{equation}\label{06-08-2}
g(s)\le c_6s^{q_0}+\epsilon,\quad \forall s\ge0,
\end{equation}
for some $q_0\in(0,1]$, $c_6>0$ and $\epsilon>0$.
Assume more that  $c_6$ is small enough when $q_0=1$.

Then problem (\ref{eq1.1bt}) admits a  weak nonnegative solution $u_\mu$ which satisfies
\begin{equation}\label{1.8}
 u_\mu\ge \varrho  \mathbb{G}_\alpha[ w_\mu].
\end{equation}

$(ii)$ Suppose that
\begin{equation}\label{ag10}
 g(s)\le c_7s^{q_*}+\epsilon,\quad \forall s\in[0,1]
\end{equation}
and
\begin{equation}\label{ag1+}
g_{\infty}:=\int_1^{+\infty} g(s)s^{-1-p^*_{\beta}}ds<+\infty,
\end{equation}
where $c_7, \epsilon>0$, $q_*>1$ and $p^*_{\beta}=\frac{N+\beta}{N-2\alpha+\beta}$.

Then there exist $\sigma_0,\varrho_0>0$  depending on $c_7, q_*, g_\infty$ and $ p^*_\beta$  such that for $\varrho\in[0,\varrho_0)$ and $\epsilon\in[0,\epsilon_0)$, problem
(\ref{eq1.1bt}) admits a nonnegative  weak solution $u_\mu$ which satisfies (\ref{1.8}).
\end{teo}

The rest of this paper is organized as follows. In section \S2, we recall some basic results
on Green kernel and Poisson kernel related to  the Marcinkiewicz
space. Section \S3 is addressed to prove the existence of weak solution to elliptic equation
with small forcing measure. Finally, we obtain weak solution to elliptic equation with small
 boundary type measure.

\setcounter{equation}{0}
\section{Preliminary}

In order to obtain the weak solution of (\ref{eq1.1})-(\ref{eq1.01}) with integral subcritical nonlinearity,
we have to introduce the Marcinkiewicz
space and recall some related  estimate.

\begin{definition}
Let $\Theta\subset \R^N$ be a domain and $\varpi$ be a positive
Borel measure in $\Theta$. For $\kappa>1$,
$\kappa'=\kappa/(\kappa-1)$ and $u\in L^1_{loc}(\Theta,d\mu)$, we
set
\begin{equation}\label{mod M}
\|u\|_{M^\kappa(\Theta,d\varpi)}=\inf\left\{c\in[0,\infty]:\int_E|u|d\varpi\le
c\left(\int_Ed\varpi\right)^{\frac1{\kappa'}},\ \forall E\subset \Theta,\,E\
{\rm Borel}\right\}
\end{equation}
and
\begin{equation}\label{spa M}
M^\kappa(\Theta,d\varpi)=\{u\in
L_{loc}^1(\Theta,d\varpi):\|u\|_{M^\kappa(\Theta,d\varpi)}<+\infty\}.
\end{equation}
\end{definition}

$M^\kappa(\Theta,d\varpi)$ is called the Marcinkiewicz space of
exponent $\kappa$, or weak $L^\kappa$-space and
$\|.\|_{M^\kappa(\Theta,d\varpi)}$ is a quasi-norm. We observe that
\begin{equation}\label{p m0}
\|u+v\|_{M^\kappa(\Theta,d\varpi)}\le\|u\|_{M^\kappa(\Theta,d\varpi)}+\|v\|_{M^\kappa(\Theta,d\varpi)}
\end{equation}
and
\begin{equation}\label{p m1}
\|t u\|_{M^\kappa(\Theta,d\varpi)}=t\|u\|_{M^\kappa(\Theta,d\varpi)},\quad \forall t>0.
\end{equation}

\begin{proposition}\label{pr 1} \cite{BBC,CC}
Assume that $1\le q< \kappa<\infty$ and $u\in L^1_{loc}(\Theta,d\varpi)$.
Then there exists  $c_8>0$ dependent of $q,\kappa$ such that
$$\int_E |u|^q d\varpi\le c_8\|u\|_{M^\kappa(\Theta,d\varpi)}\left(\int_E d\varpi\right)^{1-q/\kappa}$$
for any Borel set $E$ of $\Theta$.
\end{proposition}

The next estimate is the key-stone in the proof of Theorem \ref{teo 1} to control the nonlinearity in $\{g\ge 1\}$.
\begin{proposition}\label{general}
Let $\alpha\in(0,1]$, $\beta\in[0,\alpha]$ and $p_\beta^*=\frac{N+\beta}{N-2\alpha+\beta}$, then there exists $c_9>0$ such that
\begin{equation}\label{annex 0}
\|\mathbb{G}_\alpha[\nu]\|_{M^{p_\beta^*}(\Omega,\rho^\beta
dx)}\le c_9\|\nu\|_{\mathfrak{M}(\Omega,\rho^\beta)}.
\end{equation}

\end{proposition}
{\bf Proof.} When $\alpha\in(0,1)$, it follows by
\cite[Proposition 2.2]{CV1} that for $\gamma\in[0,\alpha]$, there exists $c_{10}>0$ such that
$$
\|\mathbb{G}_\alpha[\nu]\|_{M^{k_{\alpha,\beta,\gamma}}(\Omega,\rho^\gamma
dx)}\le c_{10}\|\nu\|_{\mathfrak{M}(\Omega,\rho^\beta)},
$$
where
$$
k_{\alpha,\beta,\gamma}= \left\{ \arraycolsep=1pt
\begin{array}{lll}
\frac{N+\gamma}{N-2\alpha+\beta},\quad &{\rm if}\ \gamma<
\frac{N\beta}{N-2\alpha},\\[2mm]
\frac{N}{N-2\alpha},\quad &{\rm if\ not}.
\end{array}
\right.
$$
We just take $\gamma=\beta$, then $k_{\alpha,\beta,\gamma}=p_\beta^*$ and (\ref{annex 0}) holds.

  When $\alpha=1$, (\ref{annex 0}) follows by \cite[Theorem 3.5]{V}.
 \qquad$\Box$

\smallskip
The following proposition does not just provide regularity but also plays an essential role to control in $\{g<1\}$.

\begin{proposition}\label{pr5}   Let $\alpha\in(0,1]$ and $\beta\in [0, \alpha]$, then the mapping $f\mapsto \mathbb G_\alpha[f]$ is compact from  $L^{1}(\Omega,\rho^\beta dx)$ into $L^{q}(\Omega)$ for any $q\in [1,\frac{N}{N+\beta-2\alpha})$.
Moreover, for  $q\in [1,\frac{N}{N+\beta-2\alpha})$, there exists $c_{11}>0$ such that for any $f\in L^{1}(\Omega,\rho^{\beta}dx)$
  \begin{equation}\label{power1}
  \norm{\mathbb G_\alpha[f]}_{L^q(\Omega)}\leq c_{11}\norm f_{L^{1}(\Omega,\rho^{\beta}dx)}.
\end{equation}
 \end{proposition}
{\bf Proof.} When $\alpha\in(0,1)$ and  $\beta\in [0, \alpha]$, it follows by \cite[Proposition 2.5]{CV1} that  for $p\in (1,\frac{N}{N-2\alpha+\beta})$, there exists $c_{12}>0$ such that for any $f\in L^{1}(\Omega,\rho^{\beta}dx)$
  \begin{equation}\label{power2}
  \norm{\mathbb G_\alpha[f]}_{W^{2\alpha-\gamma,p}(\Omega)}\leq c_{12}\norm f_{L^{1}(\Omega,\rho^{\beta}dx)},
\end{equation}
where $\gamma=\beta+\frac{N(p-1)}{p}$ if $\beta>0$ and
$\gamma>\frac{N(p-1)}{p}$ if $\beta=0$.
By \cite[Theorem 6.5]{NPV}, the embedding of $W^{2\alpha-\gamma,p}(\Omega)$ into $L^{q}(\Omega)$ is compact, then
the mapping $f\mapsto \mathbb G_\alpha[f]$ is compact from  $L^{1}(\Omega,\rho^\beta dx)$ into $L^{q}(\Omega)$ for any $q\in [1,\frac{N}{N+\beta-2\alpha})$.
We observe that (\ref{power1}) follows by (\ref{power2}) and the embedding inequality.

When $\alpha=1$ and $\beta\in[0,1]$, it  follows by   \cite[Theorem 2.7]{BV} that
\begin{equation}\label{power3}
  \norm{\mathbb G_\alpha[f]}_{W_0^{1,\frac{N}{N-1+\beta}}(\Omega)}\leq c_{13}\norm f_{L^{1}(\Omega,\rho^{\beta}dx)},
\end{equation}
where $c_{13}>0$.
 By the compactness of the embedding from $W_0^{1,\frac{N}{N-1+\beta}}(\Omega)$ into $L^{q}(\Omega)$ with $q\in [1,\frac{N}{N+\beta-2})$, we have that
the mapping $f\mapsto \mathbb G_\alpha[f]$ is compact from  $L^{1}(\Omega,\rho^\beta dx)$ into $L^{q}(\Omega)$ for  $q\in [1,\frac{N}{N+\beta-2})$.
Similarly, (\ref{power1}) follows by (\ref{power3}) and the related embedding inequality.
\qquad$\Box$

\smallskip

When we deal with problem (\ref{eq1.1b}), the Poisson kernel changes the boundary measure to forcing term and the following proposition
plays an important role in obtaining the weak solution to (\ref{1.7}) replaced $g_n$ by $g$.
\begin{proposition}\cite[Theorem 2.5]{BV}\label{pr 2}
Let $\gamma>-1$ and  $p_\gamma=\frac{N+\gamma}{N-1}$, then there exists $c_{14}>0$ such that
\begin{equation}
\|\mathbb{P}[\nu]\|_{M^{p_\gamma}(\Omega,\rho^\gamma
dx)}\le c_{14}\|\nu\|_{\mathfrak{M}^b(\partial\Omega)}.
\end{equation}

\end{proposition}

\setcounter{equation}{0}
\section{Forcing measure}

\subsection{ Sub-linear }
In this subsection, we are devoted to prove the existence of weak solution to (\ref{eq1.1})
when the nonlinearity is sub-linear.
\smallskip

\noindent{\bf Proof of Theorem \ref{teo 1} part $(i)$.}
Let $\beta\in[0,\alpha]$, we define the space
$$C_{\beta}(\bar \Omega)=\{\zeta\in C(\bar \Omega):\rho^{-\beta}\zeta\in C(\bar \Omega)\}$$
endowed with the norm
$$\norm{\zeta}_{C_{\beta}(\bar\Omega)}=\|\rho^{-\beta}\zeta\|_{C(\bar\Omega)}. $$
Let  $\{\nu_n\}\subset C^1(\bar \Omega)$ be a sequence of nonnegative functions such that
$\nu_{n }\to\nu $ in  sense of duality with $C_{\beta}(\bar
\Omega)$, that is,
\begin{equation}\label{06-08}
  \lim_{n\to\infty}\int_{\bar \Omega}\zeta \nu_{n }dx=\int_{\bar \Omega}\zeta d\nu,\qquad\forall \zeta\in C_{\beta}(\bar \Omega).
\end{equation}
By the
Banach-Steinhaus Theorem, $\norm{\nu_{n}}_{\mathfrak M
(\Omega,\rho^\beta)}$ is bounded independently of $n$. We may assume that $\norm{\nu_n}_{L^1(\Omega,\rho^{\beta}dx)}\le \norm{\nu}_{\mathfrak M (\Omega,\rho^{\beta})}=1$ for all $n\ge1$.  We consider a sequence $\{g_n\}$ of $C^1$ nonnegative  functions defined on $\R_+$
such that $g_n(0)=g(0)$,
\begin{equation}\label{06-08-1}
  g_n\le g_{n+1}\le g,\quad \sup_{s\in\R_+}g_n(s)=n\quad{\rm and}\quad \lim_{n\to\infty}\norm{g_n-g}_{L^\infty_{loc}(\R_+)}=0.
\end{equation}
We set
$$M(v)=\norm{v}_{L^{1}(\Omega)}.$$

{\it Step 1. To prove that for $n\geq 1$,
 \begin{equation}\label{002.3}
 \arraycolsep=1pt
\begin{array}{lll}
 (-\Delta)^\alpha u= g_{n}(u)+\sigma\nu_n\quad & {\rm in}\quad\Omega,\\[2mm]
 \phantom{   (-\Delta)^\alpha }
u=0\quad & {\rm in}\quad \Omega^c
\end{array}
 \end{equation}
admits a nonnegative solution $u_n$ such that
$$M(u_n)\le \bar\lambda,$$
where $\bar\lambda>0$ independent of $n$. }

To this end, we define the operators $\{\mathcal{T}_n\}$ by
 $$\mathcal{T}_nu=\mathbb{G}_\alpha\left[g_n(u)+\sigma \nu_n\right],\qquad \forall u\in L^1_+(\Omega),$$
 where $L^1_+(\Omega)$ is the positive cone of $L^1(\Omega)$.
By  (\ref{power1})  and (\ref{06-08-2}), we have that
\begin{equation}\label{23-05-0}
 \arraycolsep=1pt
\begin{array}{lll}
  M(\mathcal{T}_nu)\le c_{11}\norm{g_n(u)+\sigma \nu_n}_{L^1 (\Omega,\rho^{\beta}dx)}
  \\[2mm] \phantom{---- }
  \le c_2c_{11} \int_{\Omega}u^{p_0}\rho^\beta(x)dx+c_6(\sigma+\epsilon)
  \\[2mm] \phantom{---- }
  \le c_2c_{15}\int_{\Omega}u^{p_0} dx+c_6(\sigma+\epsilon)
  \\[2mm] \phantom{---- }
  \le c_2 c_{16}(\int_{\Omega}u dx)^{p_0}+c_6(\sigma+\epsilon)
  \\[2mm] \phantom{---- }
=c_2c_{16}M(u)^{p_0}+c_6(\sigma+\epsilon),
\end{array}
\end{equation}
where $c_{15},c_{16}>0$ independent of $n$.
Therefore, we derive that
$$
  M(\mathcal{T}_nu)\le c_2c_{16} M(u)^{p_0}+c_{11}(\sigma+\epsilon).
$$

If we assume that $M(u)\le \lambda$ for some $\lambda>0$, it implies
$$
 M(\mathcal{T}_nu)\le c_2c_{16} \lambda^{p_0}+c_{11}(\sigma+\epsilon).
$$
 In the case of $p_0<1$,   the equation
$$
c_2c_{16}\lambda^{p_0 }+c_{11}(\sigma+\epsilon)=\lambda
$$
admits a unique positive root $\bar\lambda$.
In the case of $p_0=1$,  for $c_2>0$ satisfying $c_2c_{16}<1$, the equation
$$
c_2c_{16}\lambda+c_{11}(\sigma+\epsilon)=\lambda
$$
admits a unique positive root $\bar\lambda$.
For $M(u)\le \bar\lambda$, we obtain that
  \begin{equation}\label{07-05-5jingxuan}
  M(\mathcal{T}_nu)\le c_2c_{16}\bar\lambda^{p_0}+c_{11}(\sigma+\epsilon)= \bar\lambda.
  \end{equation}
Thus, $\mathcal{T}_n$ maps $L^1(\Omega)$ into itself. Clearly, if $u_m\to u$ in $L^1(\Omega)$ as $m\to\infty$, then $g_n(u_m)\to g_n(u)$ in $L^1(\Omega)$ as $m\to\infty$, thus $\mathcal{T}_n$ is continuous.
For any fixed $n\in\N$, $\mathcal{T}_nu_m=\mathbb{G}_\alpha\left[g_n(u_m)+\sigma \nu_n\right]$ and $\{g_n(u_m)+\sigma \nu_n\}_m$ is
uniformly bounded in $L^1(\Omega,\rho^\beta dx)$, then it follows by Proposition \ref{pr5} that
$\{\mathbb{G}_\alpha\left[g_n(u_m)+\sigma \nu_n\right]\}_m$ is pre-compact in $L^1(\Omega)$, which implies that
$\mathcal{T}_n$ is a compact operator.

Let
$$
\displaystyle\begin{array}{lll}\displaystyle
\mathcal{G}=\{u\in L^1_+(\Omega): \ M(u)\le \bar\lambda \},
\end{array}
$$
  which is a closed and convex
set of $L^1(\Omega)$.  It infers by (\ref{07-05-5jingxuan}) that
$$\mathcal{T}_n(\mathcal{G})\subset \mathcal{G}.$$
 It follows by Schauder's fixed point theorem that there exists some $u_n\in L^1_+(\Omega)$ such that
$\mathcal{T}_nu_n=u_n$ and $M(u_n)\le \bar\lambda,$
where $\bar\lambda>0$ independent of $n$.

We observe that $u_n$ is a classical solution of  (\ref{002.3}). For $\alpha=1$, since $g_n$ bounded and $C^1$, then it is natural to see that. When $\alpha\in(0,1)$,  let  open set $O$ satisfy $ O\subset \bar O\subset \Omega$.
By  \cite[Proposition 2.3]{RS}, for $\theta\in(0,2\alpha)$, there exists $c_{17}>0$ such that
$$\norm{u_n}_{C^{\theta}(O)}\le c_{17}\{\norm{g(u_n)}_{L^\infty(\Omega)}+\sigma\norm{\nu_n}_{L^{\infty}(\Omega)}\},$$
then applied  \cite[Corollary 2.4]{RS}, $u_n$ is $C^{2\alpha+\epsilon_0}$ locally in $\Omega$ for some $\epsilon_0>0$.
Then $u_n$ is a classical solution of (\ref{002.3}).
Moreover, from \cite[Lemma 2.2]{CV2}, we derive that
\begin{equation}\label{5.60000}
\int_\Omega u_n(-\Delta)^\alpha\xi dx=\int_\Omega g_n(u_n)\xi dx+\sigma\int_\Omega\xi
\nu_ndx,\quad \forall\xi\in \mathbb{X}_{\alpha}.
\end{equation}

%%%%%%%%%%%%%%%%%%%%%%%%%%%%%%%%%%%%%%%%%%%%%%%%%%%%%%%%%%%%%%%%%%%%%%%%%%%%%%%%%%%%%%%%%%%%%%%%%%%%%%%%%%%%%%%%%%%%%%%%%%%%%%%%%%%%%%%%%%%%%%%%%%%%%%%%%%%%%%%

{\it Step 2. Convergence. } We observe that  $\{g_n( u_n)\}$ is uniformly bounded in $L^1(\Omega,\rho^\beta dx)$, so is $\{\nu_n\}$.
By Proposition \ref{pr5}, there exist a subsequence $\{u_{n_k}\}$ and $u$ such that
$u_{n_k}\to u$ a.e. in $\Omega$ and in $L^1(\Omega)$, then by (\ref{06-08-2}), we derive that
$g_{n_k}(u_{n_k}) \to g( u)$  in $L^1(\Omega)$. Pass the limit of (\ref{5.60000}) as $n_k\to \infty$ to derive that
 $$\int_\Omega u(-\Delta)^\alpha\xi=\int_\Omega g(u)\xi dx+\sigma\int_\Omega\xi d\nu,\quad \forall \xi\in\mathbb{X}_\alpha, $$
thus $u$ is a weak solution of (\ref{eq1.1})-(\ref{eq1.01}) and $u$ is nonnegative since $\{u_n\}$ are nonnegative.
\qquad$\Box$

\subsection{Integral subcritical }
In this subsection, we prove   the existence of weak solution to (\ref{eq1.1})
when the nonlinearity is integral subcritical. We first introduce an auxiliary lemma.
\begin{lemma}\label{lm 08-09}
Assume that $g:\R_+\mapsto\R_+$ is a continuous function satisfying
\begin{equation}\label{p}
\int_1^{+\infty} g(s)s^{-1-p}ds<+\infty
\end{equation}
for some $p>0$.
 Then there is a sequence
real positive numbers $\{T_n\}$ such that
$$\lim_{n\to\infty}T_n=\infty\quad{\rm and}\quad \lim_{n\to\infty}g(T_n)T_n^{-p}=0.$$

\end{lemma}
{\it Proof.} Let $\{s_n\}$ be a sequence of real positive numbers converging to $\infty$. We observe
\begin{eqnarray*}
\int_{s_n}^{2s_n}g(t)t^{-1-p}dt&\ge&
\min_{t\in[s_n,2s_n]}g(t)(2s_n)^{-1-p}\int_{s_n}^{2s_n}dt
\\&=&2^{-1-p}s_n^{-p}\min_{t\in[s_n,2s_n]}g(t)
\end{eqnarray*}
and by (\ref{p}),
\begin{eqnarray*}
\lim_{n\to\infty}\int_{s_n}^{2s_n}g(t)t^{-1-p}dt=0.
\end{eqnarray*}
Then we choose $T_n\in[s_n,2s_n]$ such that $g(T_n)=\min_{t\in[s_n,2s_n]}g(t)$ and then
the claim follows.
\qquad$\Box$

\smallskip

\noindent{\bf Proof of Theorem \ref{teo 1} part $(ii)$.}
Let  $\{\nu_n\}\subset C^1(\bar \Omega)$ be a sequence of nonnegative functions such that
$\nu_{n }\to\nu $ in  sense of duality with $C_{\beta}(\bar
\Omega)$ and we may assume that $\norm{\nu_n}_{L^1(\Omega,\rho^{\beta}dx)}\le 2\norm{\nu}_{\mathfrak M (\Omega,\rho^{\beta})}=1$ for all $n\ge1$.
We consider a sequence $\{g_n\}$ of $C^1$ nonnegative  functions defined on $\R_+$
satisfying $g_n(0)=g(0)$ and (\ref{06-08-1}).
We set
$$M_1(v)=\norm{v}_{M^{p_\beta^*}(\Omega,\rho^\beta dx)}\quad{\rm and}\quad M_2(v)=\norm{v}_{L^{p_*}(\Omega)},$$
where $p_\beta^*$ and $p_*$ are from (\ref{1.1}) and (\ref{1.4}).
We may assume that $p_*\in(1, \frac{N}{N-2\alpha+\beta})$.
In fact, if $p_*\ge \frac{N}{N-2\alpha+\beta}$, then  for any given $p\in(1, \frac{N}{N-2\alpha+\beta})$, (\ref{1.4}) implies that
$$g(s)\le c_3s^p+\epsilon,\quad \forall s\in[0,1].$$

{\it Step 1. To prove that for $n\geq 1$,
 \begin{equation}\label{2.3}
 \arraycolsep=1pt
\begin{array}{lll}
 (-\Delta)^\alpha u= g_{n}(u)+\sigma\nu_n\quad & {\rm in}\quad\Omega,\\[2mm]
 \phantom{   (-\Delta)^\alpha_n }
u=0\quad & {\rm in}\quad \Omega^c
\end{array}
 \end{equation}
admits a nonnegative solution $u_n$ such that
$$M_1(u_n)+M_2(u_n)\le \bar\lambda,$$
where $\bar\lambda>0$ independent of $n$. }

To this end, we define the operators $\{\mathcal{T}_n\}$ by
 $$\mathcal{T}_nu=\mathbb{G}_\alpha\left[g_n(u)+\sigma \nu_n\right],\qquad \forall u\in L^1_+(\Omega).$$
By  Proposition \ref{general}, we have
\begin{eqnarray}
  M_1(\mathcal{T}_nu) &\le& c_9\norm{g_n(u)+\sigma \nu_n}_{L^1 (\Omega,\rho^{\beta}dx)}\nonumber\\[2.5mm]
   &\le & c_9 [\norm{g_n(u)}_{L^1(\Omega,\rho^{\beta}dx)}+\sigma].\label{06-08-10}
\end{eqnarray}
In order to deal with $\norm{g_n(u)}_{L^1(\Omega,\rho^{\beta}dx)}$, for $\lambda
>0$ we set $S_\lambda=\{x\in\Omega:u(x)>\lambda\}$  and
$\omega(\lambda)=\int_{S_\lambda}\rho^\beta dx$,
\begin{equation}\label{chenyuhang1}
\displaystyle\begin{array}{lll}
\displaystyle\norm{g_n(u)}_{L^1(\Omega,\rho^{\beta}dx)}\le \int_{S^c_1}g(u)\rho^{\beta}dx+\int_{S_1}
g(u)\rho^{\beta}dx.
%\\[4mm]\phantom{-------\ }
%\displaystyle\leq  c_2\int_{S_1^c}u^{p_*}\rho^{\beta}dx+\int_{S_1} g(u)\rho^{\beta}dx.
%\leq  g(1)\int_{\Omega}\rho^{\beta}dx+\omega(1) g(1)+\int_1^\infty \omega(s)dg(s).
\end{array}
\end{equation}
We first deal with $\int_{S_1} g(u)\rho^{\beta}dx$.  In fact, we observe that
$$\int_{S_1} g(u)\rho^{\beta}dx=\omega(1) g(1)+\int_1^\infty \omega(s)dg(s),$$
where
$$\int_1^\infty  g(s)d\omega(s)=\lim_{T\to\infty}\int_1^T   g(s)d\omega(s).
$$
It infers by  Proposition \ref{pr 1} and Proposition \ref{general} that
there exists  $c_{18}>0$ such that
\begin{equation}\label{2.4}
\omega(s)\leq c_{18}M_1(u)^{p_{\beta}^*}s^{-p_{\beta}^*}
\end{equation}
and by  (\ref{1.4}) and Lemma \ref{lm 08-09} with $p=p^*_\beta$, there exist a sequence of increasing numbers $\{T_j\}$ such that $T_1>1$ and
$T_j^{-p_\beta^*}  g(T_j)\to 0$ when $j\to\infty$, thus
$$\displaystyle\begin{array}{lll}
\displaystyle \omega(1)  g(1)+ \int_1^{T_j} \omega(s)d g(s) \le c_{18}M_1(u)^{p_{\beta}^*} g(1)+c_{18}M(u)^{p_{\beta}^*}\int_1^{T_j} s^{-p_\beta^*}d g(s)
\\[4mm]\phantom{-----}\displaystyle
\leq c_{18}M_1(u)^{p_{\beta}^*}{T_j}^{-p_\beta^*}
g(T_j)+\frac{c_{18}M_1(u)^{p_{\beta}^*}}{p_\beta^*+1}\int_1^{T_j}
s^{-1-p_\beta^*} g(s)ds.
\end{array}$$
Therefore,
\begin{equation}\label{06-08-11}
\displaystyle\begin{array}{lll}
\int_{S_1} g(u)\rho^{\beta}dx=\omega(1)g(1)+ \int_1^\infty \omega(s)\ dg(s)
\\[3mm]\phantom{------}
\leq \frac{c_{18}M_1(u)^{p_{\beta}^*}}{p_\beta^*+1}\int_1^\infty s^{-1-p_\beta^*} g(s)ds
\\[3mm]\phantom{------}
\displaystyle =c_{18}g_\infty M_1(u)^{p_{\beta}^*},
\end{array}
\end{equation}
where $c_{18}>0$ independent of $n$.

We next deal with $  \int_{S^c_1}g(u)\rho^{\beta}dx$.
For $p_*\in(1, \frac{N}{N-2\alpha+\beta})$, we have that
\begin{equation}\label{4.1}
\displaystyle\begin{array}{lll}
 \int_{S^c_1}g(u)\rho^{\beta}dx\le c_3\int_{S_1^c}u^{p_*}\rho^{\beta}dx+\epsilon\int_{S_1^c}\rho^{\beta}dx
 \\[3mm]\phantom{------}
 \le c_3c_{19}\int_{\Omega}u^{p_*}dx+c_{19}\epsilon
  \\[3mm]\phantom{------}
 \leq c_3c_{19}M_2(u)^{p_*} +c_{19}\epsilon,
\end{array}
\end{equation}
where $c_{19}>0$ independent of $n$.

Along with (\ref{06-08-10}), (\ref{chenyuhang1}), (\ref{06-08-11}) and (\ref{4.1}), we derive
\begin{equation}\label{05-09-4}
  M_1(\mathcal{T}_nu)\le c_9c_{18}g_\infty M_1(u)^{p_{\beta}^*}+c_9c_3c_{19}M_2(u)^{p_*}+c_9c_{19}\epsilon+c_9\sigma.
\end{equation}
By \cite[Theorem 6.5]{NPV} and (\ref{power1}), we derive that
$$
\displaystyle\begin{array}{lll}
M_2(\mathcal{T}_nu)\le c_{11}\norm{g_n(u)+\sigma \nu_n}_{L^1 (\Omega,\rho^{\beta}dx)},
\end{array}
$$
which along with  (\ref{chenyuhang1}), (\ref{06-08-11}) and (\ref{4.1}), implies that
\begin{equation}\label{4.3}
  M_2(\mathcal{T}_nu)\le c_{11}c_{18}g_\infty M_1(u)^{p_{\beta}^*}+c_{11}c_3c_{19}M_2(u)^{p_*}+c_{11}c_{19}\epsilon+c_{11}\sigma.
\end{equation}

Therefore,  inequality (\ref{05-09-4}) and (\ref{4.3}) imply that
$$
M_1(\mathcal{T}_nu)+M_2(\mathcal{T}_nu)\le c_{20}g_\infty M_1(u)^{p_{\beta}^*}+c_{13}M_2(u)^{p_*}+c_{21}\epsilon+c_{22}\sigma,
$$
where $c_{20}=(c_9+c_{11})c_{18}$, $c_{21}=(c_9+c_{11})c_{19}$ and $c_{22}=c_9+c_{11}$.
If we assume that $M_1(u)+M_2(u)\le \lambda$, implies
$$
 M_1(\mathcal{T}_nu)+M_2(\mathcal{T}_nu)\le c_{20}g_\infty\lambda^{p_{\beta}^*}+c_{21}\lambda^{p_*}+c_{21}\epsilon+c_{22}\sigma.
$$
Since $p_{\beta}^*,\ p_*>1$, then there exist $\sigma_0>0$ and $\epsilon_0>0$ such that for any $\sigma\in(0,\sigma_0]$ and $\epsilon\in(0,\epsilon_0]$,  the equation
$$
c_{20}g_\infty\lambda^{p_{\beta}^*}+c_{21}\lambda^{p_*}+c_{21}c_2\epsilon+c_{22}\sigma=\lambda
$$
admits the largest root $\bar\lambda>0$.

We redefine $M(u)=M_1(u)+M_2(u)$, then for  $M(u)\le \bar\lambda$, we obtain that
 \begin{equation}\label{2.2}
  M(\mathcal{T}_nu)\le c_{20}g_\infty\bar\lambda^{p_{\beta}^*}+c_{21}\bar\lambda^{p_*}+c_{21}\epsilon+c_{22}\sigma= \bar\lambda.
 \end{equation}
Especially, we have that
$$\norm{\mathcal{T}_nu}_{L^1(\Omega)}\le c_8M_1(\mathcal{T}_nu)|\Omega|^{1-\frac1{p_\beta^*}}\le c_{23} \bar\lambda\quad{\rm if}\quad M(u)\le \bar\lambda.$$
Thus, $\mathcal{T}_n$ maps $L^1(\Omega)$ into itself. Clearly, if $u_m\to u$ in $L^1(\Omega)$ as $m\to\infty$, then $g_n(u_m)\to g_n(u)$ in $L^1(\Omega)$ as $m\to\infty$, thus $\mathcal{T}_n$ is continuous.
For any fixed $n\in\N$, $\mathcal{T}_nu_m=\mathbb{G}_\alpha\left[g_n(u_m)+\sigma \nu_n\right]$ and $\{g_n(u_m)+\sigma \nu_n\}_m$ is
uniformly bounded in $L^1(\Omega,\rho^\beta dx)$, then it follows by Proposition \ref{pr5} that
$\{\mathbb{G}_\alpha\left[g_n(u_m)+\sigma \nu_n\right]\}_m$ is pre-compact in $L^1(\Omega)$, which implies that
$\mathcal{T}_n$ is a compact operator.

Let
$$
\displaystyle\begin{array}{lll}\displaystyle
\mathcal{G}=\{u\in L^1_+(\Omega): \ M(u)\le \bar\lambda \}
\end{array}
$$
  which is a closed and convex
set of $L^1(\Omega)$.  It infers by (\ref{2.2}) that
$$\mathcal{T}_n(\mathcal{G})\subset \mathcal{G}.$$
 It follows by Schauder's fixed point theorem that there exists some $u_n\in L^1_+(\Omega)$ such that
$\mathcal{T}_nu_n=u_n$ and $M(u_n)\le \bar\lambda,$
where $\bar\lambda>0$ independent of $n$.

In fact, $u_n$ is a classical solution of  (\ref{2.3}).  For $\alpha=1$, since $g_n$ bounded and $C^1$, then it is natural to see that. When $\alpha\in(0,1)$, let  open set $O$ satisfy $ O\subset \bar O\subset \Omega$.
By  \cite[Proposition 2.3]{RS}, for $\theta\in(0,2\alpha)$, there exists $c_{24}>0$ such that
$$\norm{u_n}_{C^{\theta}(O)}\le c_{24}\{\norm{g(u_n)}_{L^\infty(\Omega)}+\sigma\norm{\nu_n}_{L^{\infty}(\Omega)}\},$$
then applied  \cite[Corollary 2.4]{RS}, $u_n$ is $C^{2\alpha+\epsilon_0}$ locally in $\Omega$ for some $\epsilon_0>0$.
Then $u_n$ is a classical solution of (\ref{2.3}).
Moreover,
\begin{equation}\label{5.6}
\int_\Omega u_n(-\Delta)^\alpha\xi dx=\int_\Omega g_n(u_n)\xi dx+\sigma\int_\Omega\xi
\nu_ndx,\quad \forall\xi\in \mathbb{X}_{\alpha}.
\end{equation}
%%%%%%%%%%%%%%%%%%%%%%%%%%%%%%%%%%%%%%%%%%%%%%%%%%%%%%%%%%%%%%%%%%%%%%%%%%%%%%%%%%%%%%%%%%%%%%%%%%%%%%%%%%%%%%%%%%%%%%%%%%%%%%%%%%%%%%%%%%%%%%%%%%%%%%%%%%%%%%%

{\it Step 2. Convergence. } Since  $\{g_n( u_n)\}$ and $\{\nu_n\}$ are uniformly bounded in $L^1(\Omega,\rho^\beta dx)$,
then by  Propostion \ref{pr5},  there exist a subsequence $\{u_{n_k}\}$ and $u$ such that
$u_{n_k}\to u$ a.e. in $\Omega$ and in $L^1(\Omega)$, and
$g_{n_k}(u_{n_k}) \to g( u)$ a.e. in $\Omega$.

 Finally we prove that  $g_{n_k}( u_{n_k})\to g( u)$ in $L^1(\Omega,\rho^\beta dx)$.
For $\lambda
>0$, we set $S_\lambda=\{x\in\Omega:|u_{n_k}(x)|>\lambda\}$  and
$\omega(\lambda)=\int_{S_\lambda}\rho^{\beta}dx$, then for any Borel
set $E\subset\Omega$, we have that
\begin{equation}\label{chenyuhang1000}
\displaystyle\begin{array}{lll}
\displaystyle\int_{E}|g_{n_k}(u_{n_k})|\rho^\beta dx=\int_{E\cap
S^c_{\lambda}}g(u_{n_k})\rho^\beta dx+\int_{E\cap S_{\lambda}}
g(u_{n_k})\rho^\beta dx\\[4mm]\phantom{\int_{E}|g(u_{n_k})|\rho^{\alpha}dx}
\displaystyle\leq \tilde g(\lambda)\int_E\rho^\beta dx+\int_{S_{\lambda}} g(u_{n_k})\rho^\beta dx\\[4mm]\phantom{\int_{E}g(u_{n_k})\rho^\beta dx}
\displaystyle\leq \tilde
g(\lambda)\int_E\rho^\beta dx+\omega(\lambda) g(\lambda)+\int_{\lambda}^\infty \omega(s)d
g(s),
\end{array}
\end{equation}
where $\tilde g(\lambda)=\max_{s\in[0,\lambda]}g(s)$.

On the other hand,
$$\int_{\lambda}^\infty g(s)d\omega(s)=\lim_{T_m\to\infty}\int_{\lambda}^{T_m} g(s)d\omega(s).
$$
where $\{T_m\}$ is a sequence increasing number such that
$T_m^{-p_\beta^*} g(T_m)\to 0$ as $m\to\infty$,
which could obtained by assumption (\ref{1.4}) and Lemma \ref{lm 08-09} with $p=p^*_\beta$.

It infers by  (\ref{2.4}) that
$$\displaystyle\begin{array}{lll}
\displaystyle \omega(\lambda)  g(\lambda)+ \int_{\lambda}^{T_m} \omega(s)d g(s) \le c_{18}  g(\lambda)\lambda^{-p_{\beta}^*}+c_{25}\int_{\lambda}^{T_m} s^{-p_\beta^*}d  g(s)
\\[4mm]\phantom{-----\ \int_{\lambda}^{T_m}  g(s)d\omega(s)}\displaystyle
\leq c_{25}T_m^{-p_\beta^*}g(T_m)+\frac{c_{25}}{p_\beta^*+1}\int_{\lambda}^{T_m}
s^{-1-p_\beta^*} g(s)ds,
\end{array}$$
where $c_{25}=c_{18}p_{\beta}^*$.
Pass the limit of $m\to\infty$, we have that
$$\omega(\lambda)  g(\lambda)+ \int_{\lambda}^\infty \omega(s)\ d  g(s)\leq \frac{c_{25}}{p_\beta^*+1}\int_{\lambda}^\infty s^{-1-p_\beta^*}  g(s)ds.
$$

Notice that the above quantity on the right-hand side tends to $0$
when $\lambda\to\infty$. The conclusion follows: for any
$\epsilon>0$ there exists $\lambda>0$ such that
$$\frac{c_{17}}{p_\beta^*+1}\int_{\lambda}^\infty s^{-1-p_\beta^*} g(s)ds\leq \frac{\epsilon}{2}.
$$
Since $\lambda$ is fixed, together with (\ref{chenyuhang1}), there exists $\delta>0$ such that
$$\int_E\rho^{\beta}dx\leq \delta\Longrightarrow  g(\lambda)\int_E\rho^{\beta}dx\leq\frac{\epsilon}{2}.
$$
This proves that $\{g\circ u_{n_k}\}$ is uniformly integrable in
$L^1(\Omega,\rho^\beta dx)$. Then $g\circ u_{n_k}\to g\circ u$ in
$L^1(\Omega,\rho^\beta dx)$ by Vitali convergence theorem.

Pass the limit of (\ref{5.6}) as $n_k\to \infty$ to derive that
 $$\int_\Omega u(-\Delta)^\alpha\xi=\int_\Omega g(u)\xi dx+\sigma\int_\Omega\xi d\nu,\quad \forall \xi\in\mathbb{X}_\alpha, $$
thus $u$ is a weak solution of (\ref{eq1.1})-(\ref{eq1.01}) and $u$ is nonnegative since $\{u_n\}$ are nonnegative.
\qquad$\Box$
%%%%%%%%%%%%%%%%%%%%%%%%%%%%%%%%%%%%%%%%%%%%%%%%%%%%%%%%%%%%%%%%%%%%%%%%%%%%%%%%%%%%%%%%%%%%%%%%%%%%%%%%%%
%%%%%%%%%%%%%%%%%%%%%%%%%%%%%%%%%%%%%%%%%%%%%%%%%%%%%%%%%%%%%%%%%%%%%%%%%%%%%%%%%%%%%%%%%%%%%%%%%%%%%%%%%%
%%%%%%%%%%%%%%%%%%%%%%%%%%%%%%%%%%%%%%%%%%%%%%%%%%%%%%%%%%%%%%%%%%%%%%%%%%%%%%%%%%%%%%%%%%%%%%%%%%%%%%%%%%

\setcounter{equation}{0}
\section{Boundary type  measure}

In order to prove the elliptic problem involving boundary type measure, the idea is to change
 the boundary type measure to a forcing source.

\begin{lemma} \label{lm 41}
For $\mu\in\mathfrak{M}^b_+(\partial\Omega)$, we have that
$$\mathbb{P}[\mu] \in C^1(\Omega).$$
\end{lemma}
{\bf Proof.} It infers by \cite[Proposition 2.1]{BV} that   for $(x,y)\in\Omega\times\partial\Omega$,
$$P(x,y)\le c_N|x-y|^{1-N} \quad {\rm and}\quad  |\nabla_x P(x,y)|\le c_N|x-y|^{-N},$$
then by the formulation of $\mathbb{P}[\mu]$ we have that $\mathbb{P}[\mu] \in C^1(\Omega).$\qquad$\Box$

\begin{lemma} \label{lm 42}
 Assume that $\varrho>0$, $\mu\in \mathfrak{M}^b_+(\partial\Omega)$,  $g$ is a nonnegative function satisfying (\ref{g10}) and (\ref{g1+}), $\{g_n\}$ are a sequence of $C^1$ nonnegative  functions defined on $\R_+$
satisfying $g_n(0)=g(0)$ and (\ref{06-08-1}).

Then there exists $\varrho_0>0$ and $\epsilon_0>0$ such that for $\varrho\in[0,\varrho_0]$ and $\epsilon\in[0,\epsilon_0]$,
\begin{equation}\label{12-05-0}
 \arraycolsep=1pt
\begin{array}{lll}
 -\Delta  u= g_n(u+\varrho \mathbb{P}[\mu])\quad & {\rm in}\quad\Omega,\\[2mm]
 \phantom{   -\Delta  }
u=0\quad & {\rm on}\quad \partial\Omega
\end{array}
\end{equation}
admits a nonnegative solution $w_n$ such that
$$M_1(w_n)+M_2(w_n)\le \bar\lambda$$
for some $\bar\lambda>0$ independent of $n$,
where
$$M_1(v)=\norm{v}_{M^{q^*}(\Omega,\rho dx)}\quad{\rm and}\quad M_2(v)=\norm{v}_{L^{q_*}(\Omega)},$$
with $q_*$ and $q^*$ given in (\ref{g10}) and (\ref{g1+}) respectively.

\end{lemma}
{\bf Proof.} Without loss generality, we assume $\norm{\mu}_{\mathfrak{M}^b(\partial\Omega)}=1$ and $q_*\in(1, \frac{N+1}{N-1})$.
Redenote the operators $\{\mathcal{T}_n\}$ by
 $$\mathcal{T}_nu=\mathbb{G}_1\left[g_n(u+\varrho\mathbb{P}[\mu])\right],\qquad \forall u\in L^1_+(\Omega).$$
By  Proposition \ref{general}, we have
\begin{equation}\label{12-05-2}
\begin{array}{lll}
\displaystyle
   M_1(\mathcal{T}_nu) \le c_9\norm{g_n(u+\varrho \mathbb{P}[\mu])}_{L^1 (\Omega,\rho dx)}
   \\[4mm]\phantom{---- }
   \le c_9\norm{g(u+\varrho \mathbb{P}[\mu])}_{L^1 (\Omega,\rho dx)}
   \end{array}
\end{equation}
For $\lambda>0$, we set $S_\lambda=\{x\in\Omega:u+\varrho \mathbb{P}[\mu]>\lambda\}$  and
$\omega(\lambda)=\int_{S_\lambda}\rho  dx$,
\begin{equation}\label{12-05-3}
\displaystyle\begin{array}{lll}
\displaystyle
\norm{g(u+\varrho \mathbb{P}[\mu])}_{L^1(\Omega,\rho dx)}\le \int_{S^c_1}g(u+\varrho \mathbb{P}[\mu])\rho dx+\int_{S_1}
g(u+\varrho \mathbb{P}[\mu])\rho dx.
%\\[4mm]\phantom{-------\ }
%\displaystyle\leq  c_3\int_{S_1^c}u^{q_*}\rho^{\beta}dx+\int_{S_1} g(u)\rho^{\beta}dx.
%\leq  g(1)\int_{\Omega}\rho^{\beta}dx+\omega(1) g(1)+\int_1^\infty \omega(s)dg(s).
\end{array}
\end{equation}
We first deal with $\int_{S_1} g(u+\varrho \mathbb{P}[\mu])\rho dx$.  In fact, we observe that
$$\int_{S_1} g(u+\varrho \mathbb{P}[\mu])\rho dx=\omega(1) g(1)+\int_1^\infty \omega(s)dg(s),$$
where
$$\int_1^\infty  g(s)d\omega(s)=\lim_{T\to\infty}\int_1^T   g(s)d\omega(s).
$$
It infers by  Proposition \ref{general} and Proposition \ref{pr 2} with $\gamma=1$ that
there exists  such that
\begin{equation}\label{12-05-4}
\displaystyle\begin{array}{lll}
\omega(s)\leq c_{26}\norm{u+\varrho \mathbb{P}[\mu]}_{M^{q^*}(\Omega,\rho dx)}^{q^*}s^{-q^*}
\\[2mm]
\phantom{--\ }
\le c_{27}\left(\norm{u}_{M^{q^*}(\Omega,\rho dx)}+\norm{\varrho \mathbb{P}[\mu]}_{M^{q^*}(\Omega,\rho dx)}\right)^{q^*}s^{-q^*}\\[2mm]
\phantom{--\ }
\le c_{27} \left(M_1(u)+c_{14}\varrho \right)^{q^*}s^{-q^*}
\end{array}
\end{equation}
where $c_{26},c_{27}>0$ independent of $n$.
By  (\ref{g1+}) and Lemma \ref{lm 08-09} with $p=q^*$, there exist a sequence of increasing numbers $\{T_j\}$ such that $T_1>1$ and
$T_j^{-q^*}  g(T_j)\to 0$ when $j\to\infty$, thus
$$\displaystyle\begin{array}{lll}
\displaystyle \omega(1)  g(1)+ \int_1^{T_j} \omega(s)d g(s)
\\[4mm]\phantom{--}\displaystyle
\le c_{27}\left(M_1(u)+c_{14}\varrho \right)^{p_{\beta}^*} g(1)+c_{27}\left(M_1(u)+c_{14}\varrho \right)^{q^*}\int_1^{T_j} s^{-q^*}d g(s)
\\[4mm]\phantom{--}\displaystyle
\leq c_{27}\left(M_1(u)+c_{14}\varrho \right)^{q^*}{T_j}^{-q^*}
g(T_j)
\\[4mm]\phantom{---}+\frac{c_{27}\left(M_1(u)+c_{14}\varrho \right)^{q^*}}{q^*+1}\int_1^{T_j}
s^{-1-q^*} g(s)ds.
\end{array}$$
Therefore,
\begin{equation}\label{12-05-5}
\displaystyle\begin{array}{lll}
\int_{S_1} g(u)\rho dx=\omega(1)g(1)+ \int_1^\infty \omega(s)\ dg(s)
\\[3mm]\phantom{------}
\leq \frac{c_{27}\left(M_1(u)+c_{14}\varrho \right)^{q^*}}{q^*+1}\int_1^\infty s^{-1-q^*} g(s)ds
\\[3mm]\phantom{------}
\displaystyle \le c_{28}g_\infty  M_1(u)^{q^*}+c_{28}g_\infty \varrho^{q^*},
\end{array}
\end{equation}
where $c_{28}>0$ independent of $n$.

We next deal with $  \int_{S^c_1}g(u+\varrho \mathbb{P}[\mu])\rho dx$. For $q_*\in(1, \frac{N+1}{N-1})$, we have that
\begin{equation}\label{12-05-6}
\displaystyle\begin{array}{lll}
 \int_{S^c_1}g(u+\varrho \mathbb{P}[\mu])\rho dx\le c_5\int_{S_1^c}(u+\varrho \mathbb{P}[\mu])^{q_*}\rho dx+\epsilon\int_{S_1^c}\rho dx
 \\[3mm]\phantom{--------\ }
 \le c_5c_{29}\int_{\Omega}u^{q_*}dx+c_5c_{29}\varrho^{q_*}+c_{29}\epsilon
  \\[3mm]\phantom{--------\ }
 \leq c_5c_{29}M_2(u)^{q_*} +c_5 c_{29}\varrho^{q_*}+c_{29}\epsilon,
\end{array}
\end{equation}
where $c_{29}>0$ independent of $n$.

Along with (\ref{12-05-2}), (\ref{12-05-3}), (\ref{12-05-5}) and (\ref{12-05-6}), we derive that
\begin{equation}\label{05-09-4}
  M_1(\mathcal{T}_nu)\le c_9c_{26}g_\infty M_1(u)^{q^*}+c_9 c_5c_{29}M_2(u)^{q_*}+c_9c_{29}\epsilon+c_9l_{\varrho},
\end{equation}
where $l_{\varrho}=c_{28}g_\infty\varrho^{p^*}+c_5 c_{29}\varrho^{p_*}$.
By \cite[Theorem 6.5]{NPV} and (\ref{power1}), we derive that
$$
\displaystyle\begin{array}{lll}
M_2(\mathcal{T}_nu)\le c_{11}\norm{g(u+\varrho \mathbb{P}[\mu])}_{L^1 (\Omega,\rho dx)},
\end{array}
$$
which along with  (\ref{12-05-3}), (\ref{12-05-5}) and (\ref{12-05-6}), implies that
\begin{equation}\label{4.3}
  M_2(\mathcal{T}_nu)\le c_{11}c_{26}g_\infty M_1(u)^{q^*}+c_{11} c_5c_{29}M_2(u)^{q_*}+c_{11}c_{29}\epsilon+c_{11}l_{\varrho}.
\end{equation}

Therefore,  inequality (\ref{05-09-4}) and (\ref{4.3}) imply that
$$
M_1(\mathcal{T}_nu)+M_2(\mathcal{T}_nu)\le c_{30}g_\infty M_1(u)^{q^*}+c_{31}M_2(u)^{q_*}+c_{31}\epsilon+c_{32}l_{\varrho},
$$
where $c_{30}=(c_9+c_{11})c_{26}$, $c_{31}=(c_9+c_{11})c_{5}c_{29}$ and $c_{32}=c_9+c_{11}$.
If we assume that $M_1(u)+M_2(u)\le \lambda$, implies
$$
 M_1(\mathcal{T}_nu)+M_2(\mathcal{T}_nu)\le c_{30}g_\infty\lambda^{q^*}+c_{13}\lambda^{q_*}+c_{31}\epsilon+c_{32}l_{\varrho}.
$$
Since $q^*,q_*>1$, then there exist $\varrho_0>0$ and $\epsilon_0>0$ such that for any $\varrho\in(0,\varrho_0]$ and $\epsilon\in(0,\epsilon_0]$,  the equation
$$
c_{30}g_\infty\lambda^{q^*}+c_{31}\lambda^{q_*}+c_{31}c_5\epsilon+c_{32}l_{\varrho}=\lambda
$$
admits the largest root $\bar\lambda>0$.

We redefine $M(u)=M_1(u)+M_2(u)$, then for  $M(u)\le \bar\lambda$, we obtain that
 \begin{equation}\label{2.2}
  M(\mathcal{T}_nu)\le c_{30}g_\infty\bar\lambda^{q^*}+c_{31}\bar\lambda^{q_*}+c_{31}\epsilon+c_{32}l_{\varrho}= \bar\lambda.
 \end{equation}
Especially, we have that
$$\norm{\mathcal{T}_nu}_{L^1(\Omega)}\le c_8 M_1(\mathcal{T}_nu)|\Omega|^{1-\frac1{q^*}}\le c_{33} \bar\lambda\quad{\rm if}\quad M(u)\le \bar\lambda.$$
Thus, $\mathcal{T}_n$ maps $L^1(\Omega)$ into itself. Clearly, if $u_m\to u$ in $L^1(\Omega)$ as $m\to\infty$, then $g_n(u_m)\to g_n(u)$ in $L^1(\Omega)$ as $m\to\infty$, thus $\mathcal{T}_n$ is continuous.
For any fixed $n\in\N$, $\mathcal{T}_nu_m=\mathbb{G}_1\left[g_n(u_m+\varrho \mathbb{P}[\mu])\right]$ and $\{g_n(u_m)+\varrho \mathbb{P}[\mu]\}_m$ is
uniformly bounded in $L^1(\Omega,\rho dx)$, then it follows by Proposition \ref{pr5} that
$\{\mathbb{G}_1\left[g_n(u_m+\varrho \mathbb{P}[\mu]\right]\}_m$ is pre-compact in $L^1(\Omega)$, which implies that
$\mathcal{T}_n$ is a compact operator.

Let
$$
\displaystyle\begin{array}{lll}\displaystyle
\mathcal{G}=\{u\in L^1_+(\Omega): \ M(u)\le \bar\lambda \}
\end{array}
$$
  which is a closed and convex
set of $L^1(\Omega)$.  It infers by (\ref{2.2}) that
$$\mathcal{T}_n(\mathcal{G})\subset \mathcal{G}.$$
 It follows by Schauder's fixed point theorem that there exists some $w_n\in L^1_+(\Omega)$ such that
$\mathcal{T}_nw_n=w_n$ and $M(w_n)\le \bar\lambda,$
where $\bar\lambda>0$ does not depend on $n$.
Since $g_n$ and $ \mathbb{P}[\mu]$ are $C^1$ functions by Lemma \ref{lm 41}, then $w_n$ is a classical solution of (\ref{12-05-0})
and
$$
\int_\Omega w_n(-\Delta) \xi dx=\int_\Omega g_n(w_n+\varrho \mathbb{P}[\mu])\xi dx,\quad \forall\xi\in C_0^{1.1}(\Omega).
$$

\noindent{\bf Proof of Theorem \ref{teo 2} $(ii)$.}
It derives by Lemma \ref{lm 41} that $w_n$ is a classical solution of (\ref{12-05-0}).
Denote $u_n=w_n+\varrho \mathbb{P}[\mu]$ and then
\begin{equation}\label{7.1}
\int_\Omega u_n(-\Delta)\xi=\int_\Omega g_n(u_n)\xi dx+\varrho\int_{\partial\Omega}\frac{\partial\xi(x)}{\partial\vec n_x}d\mu(x),\quad \forall \xi\in\mathbb{X}_\alpha,
\end{equation}
 Since  $\{g_n( u_n)\}$  are uniformly bounded in $L^1(\Omega,\rho dx)$,
then by  Propostion \ref{pr5},  there exist a subsequence $\{w_{n_k}\}$ and $w$ such that
$w_{n_k}\to w$ a.e. in $\Omega$ and in $L^1(\Omega)$ and then $u_{n_k}\to u$ a.e. in $\Omega$ and in $L^1(\Omega)$
where $u=w+\varrho \mathbb{P}[\mu]$. Thus, $g_{n_k}(u_{n_k}) \to g( u)$ a.e. in $\Omega$.

Similarly to the argument in {\it Proof of Theorem \ref{teo 1} part $(ii)$ in Step 2},
we have that  $g_{n_k}( u_{n_k})\to g( u)$ in $L^1(\Omega,\rho dx)$.

Pass the limit of (\ref{7.1}) as $n_k\to \infty$ to derive that
 $$\int_\Omega u(-\Delta) \xi dx=\int_\Omega g(u)\xi dx+\varrho\int_{\partial\Omega}\frac{\partial\xi}{\partial\vec n} d\mu,\quad \forall \xi\in\mathbb{X}_\alpha, $$
thus $u$ is a weak solution of (\ref{eq1.1b}) and $u$ is nonnegative since $\{u_n\}$ are nonnegative.
\qquad$\Box$
\smallskip

\noindent{\bf Proof of Theorem \ref{teo 2} $(i)$. } It proceeds similarly to the proof of Theorem \ref{teo 1} $(i)$, so we omit here. \qquad$\Box$

\setcounter{equation}{0}
\section{Boundary type measure for $\alpha\in(0,1)$}

\subsection{Basic results}
In this subsection, we devoted to study the properties of $\mathfrak{R}_\beta$ with  $\beta\in[0,\alpha]$, see the definition \ref{5.2}.
 Here and in what follows,
we  assume that $\alpha\in(0,1)$.

\begin{lemma}
Let $1\le\beta'\le \beta\le \alpha$, then
\begin{equation}\label{5.3}
\emptyset\not=  \mathfrak{R}_{\beta'}\subset \mathfrak{R}_\beta\not=\mathfrak{M}_+(\Omega^c).
\end{equation}
\end{lemma}
{\bf Proof.} Let  $x_0\in\partial\Omega$, $x_t=x_0+t\vec{n}_{x_0}$
and $\delta_t$ be the dirac mass concentrated at $x_t$,
where $\vec{n}_{x_0}$ is the unit normal vector pointing outside at $x_0$.

Fixed $t>0$, $w_{\delta_t}(x)=|x-x_t|^{-N-2\alpha}$ for $x\in\Omega$. It is easy to
see that $w_{\delta_t}\in L^\infty(\Omega)$  and then $\delta_t\in \mathfrak{R}_\beta$ for any $\beta\in[0,\alpha]$.

 Fixed $t=0$, $w_{\delta_0}(x)=|x-x_0|^{-N-2\alpha}$ for $x\in\Omega$. We observe that $w_{\delta_0}\not\in L^1(\Omega,\rho^\alpha dx)$  and then $\delta_0\not\in \mathfrak{R}_\beta$ for any $\beta\in[0,\alpha]$.
\qquad$\Box$

\medskip

\noindent{\bf Example. }  Let  $x_0\in\partial\Omega$, $x_t=x_0+t\vec{n}_{x_0}$
and $\delta_t$ be the dirac mass concentrated at $x_t$.
 Denote
 $$\mu=\sum_{n=1}^{\infty}b_n\delta_{\frac1n},$$
  where $\{b_n\}$   a sequence nonnegative numbers will be chosen latter.
 We observe that
 $$w_\mu(x)=\sum_{n=1}^{\infty}\frac{b_n}{|x-x_{\frac1n}|^{N+2\alpha}},\quad x\in\Omega$$
 and
 $w_\mu\in L^1(\Omega,\rho^\beta dx)$ if and only if
\begin{equation}\label{5.4}
 \sum_{n=1}^{\infty} b_n n^{2\alpha-\beta}<+\infty.
\end{equation}

\begin{lemma}\label{lm 51}
Let $\mu\in \mathfrak{R}_\beta$ with $\beta\in[0,\alpha]$ and $w_\mu$ is given by (\ref{poisson type}).

$(i)$ $w_\mu\in C^1(\Omega)\cap L^1(\Omega,\rho^\beta dx)$.

$(ii)$ Let $\tilde w_\mu=\mathbb{G}_\alpha[w_\mu]$ in $\Omega$ and $\tilde w_\mu=\mu$ in $\Omega^c$, then $\tilde w_\mu$ is a weak solution of
\begin{equation}\label{5.5}
\begin{array}{lll}
(-\Delta)^\alpha u=0\quad {\rm in}\quad \Omega,
\\[2mm]\phantom{---}
u=\mu\quad {\rm in}\quad \Omega^c
\end{array}
\end{equation}
in the sense of
$$\int_\Omega u(-\Delta)^\alpha\xi dx=\int_\Omega \xi w_\mu dx,\quad\  \forall\ \xi\in \mathbb{X}_{\alpha}.$$

\end{lemma}
{\bf Proof.} $(i)$ $\mu\in \mathfrak{R}_\beta$ implies that  $w_\mu\in  L^1(\Omega,\rho^\beta dx)$ and since
the function:  $x:\to |x-y|^{-N-2\alpha}$ is $C^1(\Omega)$ for any $y\in\Omega^c$, then  $w_\mu\in C^1(\Omega)$.

 $(ii)$ For $\mu\in \mathfrak{R}_\beta$ with $\beta\in[0,\alpha]$, let $\{\mu_n\}\subset C^1_0( \R^N)$ with supp$(\mu_n)\subset \bar\Omega^c$ be a sequence of nonnegative functions such that
$\mu_{n}\to\mu $ in distribution sense.

Then we derive that
$w_{\mu_n}\in C^1(\bar\Omega)$ and there exists a unique classical solution $\mathbb{G}_\alpha[w_{\mu_n}]$ to
\begin{equation}\label{5.5}
\begin{array}{lll}
(-\Delta)^\alpha u=w_{\mu_n}\quad {\rm in}\quad \Omega,
\\[2mm]\phantom{--- \ }
u=0\qquad {\rm in}\quad \Omega^c.
\end{array}
\end{equation}
Moreover,
\begin{equation}\label{5.1.1}
\int_\Omega u(-\Delta)^\alpha\xi dx=\int_\Omega \xi w_{\mu} dx,\quad\  \forall\ \xi\in \mathbb{X}_{\alpha}.
\end{equation}
Let $u_n=\mathbb{G}_\alpha[w_{\mu_n}]+\mu_n$, then we have that
$$(-\Delta)^\alpha u_n =  (-\Delta)^\alpha \mathbb{G}_\alpha[w_{\mu_n}] + (-\Delta)^\alpha\mu_n= w_{\mu_n}-w_{\mu_n}=0$$
and (\ref{5.1.1}) holds for $u_n$. Passing the limit of $n\to\infty$, we derive that
$\tilde w_\mu$ is a weak solution of (\ref{5.5}).\qquad$\Box$

\smallskip

We note that $(i)$ Lemma \ref{lm 51}$(ii)$ indicates that $\mathbb{G}_\alpha[w_{\mu}]$ has the similar role as $\mathbb{P}[\mu]$ when $\alpha=1$;
$(ii)$ the definition \ref{weak definition bt} is equivalent to
\begin{definition}
 $u_\mu$ is a weak solution of (\ref{eq1.1bt}),
 if $u_\mu\in
L^1(\Omega)$,  $g(u_\mu)\in L^1(\Omega,\rho^\alpha dx)$  and
 $$\int_\Omega u_\mu(-\Delta)^\alpha \xi dx=\int_\Omega g(u_\mu)\xi dx+ \int_\Omega w_\mu \xi dx,\qquad \xi\in \mathbb{X}_\alpha,$$
 where $w_\mu$ is given by (\ref{poisson type}).

\end{definition}

\subsection{Proof of Theorem \ref{teo 3}}
Inspired by the proof of Theorem \ref{teo 2}, we first give an important lemma,
which is important in dealing with the  subcritical case.

\begin{lemma} \label{lm 52}
Assume that $\varrho>0$, $\mu\in \mathfrak{R}_\beta$, $g$ is a nonnegative function satisfying (\ref{ag10}) and (\ref{ag1+}),
$\{g_n\}$ are a sequence of $C^1$ nonnegative  functions defined on $\R_+$
satisfying $g_n(0)=g(0)$ and (\ref{06-08-1}).\\
Then there exists $\varrho_0>0$ and $\epsilon_0>0$ such that for $\varrho\in[0,\varrho_0]$ and $\epsilon\in[0,\epsilon_0]$,
\begin{equation}\label{21-05-0}
 \arraycolsep=1pt
\begin{array}{lll}
 (-\Delta)^\alpha  u= g_n(u+\varrho\mathbb{G}_\alpha[ w_{\mu}])\quad & {\rm in}\quad\Omega,\\[2mm]
 \phantom{ (-\Delta)^\alpha }
u=0\quad & {\rm in}\quad \Omega^c
\end{array}
\end{equation}
admits a nonnegative solution $w_n$ such that
$$M_1(w_n)+M_2(w_n)\le \bar\lambda$$
for some $\bar\lambda>0$ independent of $n$,
where
$$M_1(v)=\norm{v}_{M^{p^*_\beta}(\Omega,\rho^\beta dx)}\quad{\rm and}\quad M_2(v)=\norm{v}_{L^{q_*}(\Omega)},$$
with $q_*$ and $p^*_\beta$ given in (\ref{ag10}) and (\ref{ag1+}) respectively.

\end{lemma}
{\bf Proof.} For $\mu\in \mathfrak{R}_\beta$, we have that $w_\mu\in L^1(\Omega,\rho^\beta dx)$,
which, by Proposition \ref{pr5}, implies that $\mathbb{G}_\alpha[ w_{\mu}]\in M^{p^*_\beta}(\Omega,\rho^\beta dx)$.
It proceeds as Lemma \ref{lm 42}, replaced $\mathbb{P}[\mu]$ by $\mathbb{G}_\alpha[ w_{\mu}]$ to obtain that
 there exists $\varrho_0>0$ and $\epsilon_0>0$ such that for $\varrho\in[0,\varrho_0]$ and $\epsilon\in[0,\epsilon_0]$,
there exists $w_n$ such that
$$w_n=\mathbb{G}_\alpha[g_n(w_n+\varrho\mathbb{G}_\alpha[ w_{\mu}])]$$
and
$$M_1(w_n)+M_2(w_n)\le \bar\lambda$$
for some $\bar\lambda>0$ independent of $n$.

By Lemma \ref{lm 51} $(i)$, we see  that $w_n$ is a classical solution of (\ref{21-05-0}).
Moreover,
\begin{equation}\label{21-7.1}
\int_\Omega w_n(-\Delta)^\alpha \xi dx=\int_\Omega g_n(w_n+\varrho \mathbb{G}_\alpha[w_{\mu}])\xi dx,\quad \forall\xi\in C_0^{1.1}(\Omega).
\end{equation}

\noindent{\bf Proof of Theorem \ref{teo 3} $(ii)$.} It derives by Lemma \ref{lm 52} that $w_n$ is a classical solution of (\ref{21-05-0}).
Denote $u_n=w_n+\varrho \mathbb{G}_\alpha[w_{\mu}]$
 Since  $\{g_n( u_n)\}$  are uniformly bounded in $L^1(\Omega,\rho dx)$,
then by  Propostion \ref{pr5},  there exist a subsequence $\{w_{n_k}\}$ and $w$ such that
$w_{n_k}\to w$ a.e. in $\Omega$ and in $L^1(\Omega)$ and then $u_{n_k}\to u$ a.e. in $\Omega$ and in $L^1(\Omega)$
where $u=w+\varrho \mathbb{G}_\alpha[w_{\mu}]$. Thus, $g_{n_k}(u_{n_k}) \to g( u)$ a.e. in $\Omega$.

Similarly to the argument in {\it Proof of Theorem \ref{teo 1} part $(ii)$ in Step 2},
we have that  $g_{n_k}( u_{n_k})\to g( u)$ in $L^1(\Omega,\rho^\beta dx)$.

Pass the limit of (\ref{21-7.1}) as $n_k\to \infty$ to derive that
 $$\int_\Omega w(-\Delta)^\alpha \xi dx=\int_\Omega g(w+\varrho \mathbb{G}_\alpha[w_{\mu}])\xi dx,\quad \forall \xi\in\mathbb{X}_\alpha. $$
 Thus $u=w+\varrho \mathbb{G}_\alpha[w_{\mu}]$ is a weak solution of (\ref{eq1.1bt}) and $u$ is nonnegative since $\{w_n\}$ are nonnegative.
\qquad$\Box$
\smallskip

\noindent{\bf Proof of Theorem \ref{teo 3} $(i)$. } It proceeds similarly to the proof of Theorem \ref{teo 1} $(i)$, so we omit here. \qquad$\Box$

\end{document}